\documentclass[11pt]{amsart}
\usepackage{amssymb}
\usepackage{epsfig}
\usepackage{graphicx}
\newtheorem{theorem}{Theorem}
\newtheorem{lemma}[theorem]{Lemma}

\newtheorem{conjecture}[theorem]{Conjecture}

\newtheorem{corollary}[theorem]{Corollary}
\newtheorem{problem}[theorem]{Problem}

\title{Burnside kei}
\author{Maciej Niebrzydowski, J\'ozef~H. Przytycki}
\begin{document}
\thispagestyle{empty}
\centerline{\psfig{figure=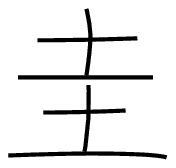,height=3cm}}

\begin{abstract}
This paper is motivated by a general question: for which values
of $k$ and $n$ is the universal Burnside kei $\bar Q(k,n)$ finite?
It is known (starting from the work of M. Takasaki (1942)) 
that $\bar Q(2,n)$ is isomorphic to the dihedral quandle $Z_n$ 
and $\bar Q(3,3)$ is isomorphic to \\$Z_3 \oplus Z_3$.
In this paper we give descriptions of $\bar Q(4,3)$ and $\bar Q(3,4)$.
We also investigate some properties of arbitrary quandles satisfying 
the universal Burnside relation $a=...a*b*...*a*b$. In particular, we prove 
that the order of a finite commutative kei is a power of 3. 
Invariants of links related to Burnside kei $\bar Q(k,n)$ are invariant 
under $n$-moves.
\end{abstract}
\maketitle

\section{Introduction}
\begin{figure}
\begin{center}
\includegraphics[height=7 cm]{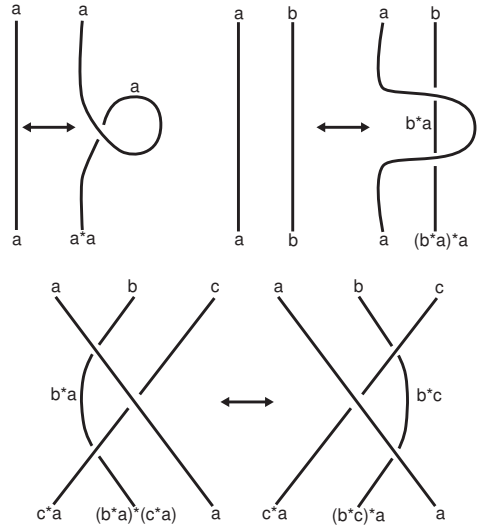}
\caption{Reidemeister moves and kei axioms.\label{rad}}
\end{center}
\end{figure}

\begin{figure}
\begin{center}
\includegraphics[height=3.5 cm]{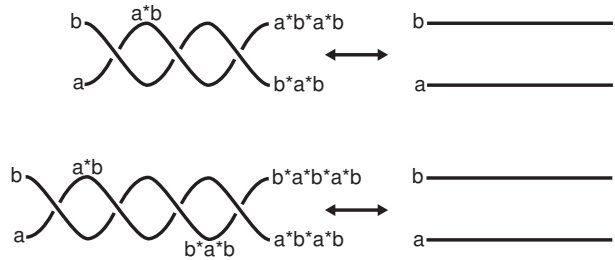}
\caption{Correspondence between $n$-moves and relation $r_n$.\label{mov}}
\end{center}
\end{figure}

Kei, \psfig{figure=dipkei.eps,height=0.3cm}, also called an involutory 
quandle, was introduced by
Mituhisa Takasaki in 1942 \cite{Tak} as 
an abstract algebra $(Q,*)$ 
with a binary operation $*: Q\times Q \to Q$ satisfying the conditions:
\begin{enumerate}
\item [(i)] $a*a=a$ for any $a\in Q$,
\item [(ii)] $(a*b)*b=a$,
\item [(iii)]$(a*b)*c= (a*c)*(b*c)$ (the right distributivity property).
\end{enumerate}
We adopt the standard convention (the left normed convention) 
 that omission of parentheses denotes the
left association, for example $a*b*c$ denotes $(a*b)*c$.
The above axioms correspond to the Reidemeister moves (see Fig.~\ref{rad}).\\
We consider free keis with the universal relation:
\begin{center}
$r_n$:$\ a=... a*b*...*a*b$,
\end{center}
in which there are $n$ letters on the right hand side and $a$, $b$ 
are any elements of the kei. We denote such kei with $k$ generators 
as $\bar Q(k,n)$ and call it the universal Burnside kei.\\
In \cite{Joy} D. Joyce associated an involutory quandle to a link. In a
similar way we can associate to every link its $n$-th Burnside kei, 
$\bar Q_n(L)$, by assigning generators to arcs of the diagram of $L$, 
writing the relation of the form $u*v=w$ for each crossing 
(here $u$ and $w$ are generators corresponding to the under-arcs 
and $v$ is assigned to the over-arc) and adding the universal 
relation $r_n$.\\
Relation $r_n$ corresponds to the local changes in the diagram 
called $n$-moves.
It follows that $\bar Q_n(L)$ is invariant under Reidemeister 
moves and $n$-moves.
For example, $\,r_3$:$\ a=b*a*b \,$ corresponds to invariance under $3$-moves 
and $\,r_4$:$\ a=a*b*a*b \,$ makes $\bar Q_4(L)$ invariant 
under $4$-moves. Fig.~\ref{mov} illustrates this correspondence 
in the case $n=3$, $4$. In fact, $\bar Q_n(L)$ is also invariant 
under rational $\frac{n}{m}$-moves \cite{D-I-P}.
\\
We notice that the relation $\,r_3$:$\ a=b*a*b \,$ is equivalent to
$\,a*b=b*a \,$, in other words, $\bar Q(k,3)$ is a free commutative 
kei on $k$ generators.\footnote{The commutative relation should 
not be confused with the abelian condition,\\
$(a*b)*(c*d)=(a*c)*(b*d)$, introduced in \cite{Joy}.}
\begin{problem}\label{1.1}
For which values of $k$ and $n$ is $\bar Q(k,n)$ a finite kei? 
How many elements does it have?
\end{problem}
In this paper we focus on finitely generated commutative keis and keis satisfying the  $4$-th universal Burnside relation $a=a*b*a*b$.
\section{Commutative keis}
\subsection{Examples}
Let us first recall that there are two well known classes of 
examples of finite commutative keis:
\begin{enumerate}
\item[(1)] dihedral kei, $Z_3$ (with $i*j=2j-i=-j-i\ modulo \ 3$), corresponding to Fox $3$-colorings, and
its direct sums $Z_3^n$ with coordinatewise operation;
\item[(2)] the third Burnside groups, 
$B(k,3)=\{x_1,...,x_k\ | w^3=1$ for any word $w \}$,
with core operation $a*b= ba^{-1}b$, and their quotients.
\end{enumerate}
Notice that $B(k,3)$ is a commutative kei as $a*b= b*a$ follows from 
$ba^{-1}b = ab^{-1}a$ which is equivalent to $(ba^{-1})^3=1$. 
Our motivation for a Problem \ref{1.1} is Burnside's theorem \cite{Bu} that 
$B(k,3)$ is a finite group.
\subsection{Some properties of commutative keis}
First let us describe some general properties of involutory quandles satisfying relation $r_3$.
Any quandle is distributive from the right, but in the case of commutative keis we also have distributivity from the left:
$$ c*(a*b)=(a*b)*c=(a*c)*(b*c)=(c*a)*(c*b).$$
From axiom (ii) in the definition of kei $Q$ it follows that:
$$\forall a,b\in Q \ \exists \ \textrm{unique}\  
c\in Q\ \textrm{such that}\  a=c*b
\ \textrm{(and obviously}\ c=a*b).$$
Here we mention that if we replace axiom (ii) with the above statement 
without the condition that $c=a*b$, we get a general definition of a quandle.
The equality $a=c*b$ is equivalent (using commutativity) to:
$$ a*c=b\ \, \textrm{and} \ \, c*a=b.$$
It follows that any commutative kei is a quasigroup\footnote{A quasigroup 
is a set $G$ together with a binary operation $\cdot$, with the property 
that for each $x$, $y\in G$, there are unique elements $w$, $z\in G$ such 
that $x\cdot w=y$ and $z\cdot x=y$.} and the set\\
$\{a,b,c=a*b\}$ is a subquandle.\\
If $m$ denotes the size of the finite commutative kei, then there 
are $ {m \choose 2} /3 $ such $3$-element subquandles and each element 
$x \in Q$ belongs to $(m-1)/2$ of them (choosing any element 
$p \in Q \setminus x$ automatically determines the third element
of the quandle, $x*p$).\\
An involutory quandle $Q$ is said to be algebraically connected 
if for each pair $a$, $b$ in $Q$, there are $a_1,a_2,\ldots ,a_s \in Q$ 
such that
$$a*a_1*a_2* \ldots *a_s=b.$$
We say that an involutory quandle is strongly algebraically connected 
if it is algebraically connected and $s=1$ in the above definition.
\begin{lemma} \label{strong}
Any kei satisfying the universal relation $r_k$, for some odd $k$, 
is strongly algebraically connected.
\end{lemma}
\begin{proof}
Our relation $r_k$ now has the form: $$a=b*a*b*\ldots *a*b.$$
Using the first axiom of quandle and the operator level relation\\
$x*y*z*y=x*(z*y)$, we can write it as:
\begin{displaymath}
a=b*b*a*b*\ldots*a*b=
\left\{ \begin{array}{l}
b*(b*a*\ldots*a*b)\\
\ \ \ \ \ \ \ \textrm{or}\\
b*(a*b*\ldots*a*b)  
\end{array} \right. , 
\end{displaymath}
depending on the length of the relation $r_k$.
In either case, in order to get from $a$ to $b$ we need to use 
only one operator (that can be written using $(k+1)/2$ letters $a$ and $b$).
\end{proof}
Every algebraically connected quandle (not necessarily involutory) is a 
metric space if we define the metric $d(x,y)$ as the minimal number 
of operators needed to obtain one element from the other. 
A significant class of algebraically connected quandles are 
knot quandles (see \cite{Joy} for a definition). Since our metric 
is unchanged under isomorphism of quandles, some metric 
properties (for example diameter) are knot invariants.
Lemma \ref{strong} states that the diameter of any quandle 
(that is, the diameter of the corresponding metric space) 
satisfying relation $r_k$ for some odd $k$ is $1$.\\
Two elements $x$ and $y$ of a quandle $Q$ are called behaviorally 
equivalent if $$z*x=z*y \ \, \textrm{for all} \ z \in Q.$$
It turns out that there are no behaviorally equivalent elements 
in quandles of the sort considered above.
\begin{lemma} \label{belements}
If $Q$ is a kei satisfying relation $r_k$, for some odd $k$, 
then it has no behaviorally equivalent elements. 
Moreover, if $z*x=z*y$ for some element $z\in Q$, then $x=y$.
\end{lemma}
\begin{proof}
Assume that $z*x=z*y$, for some $x$, $y$, $z \in Q$. 
In the case of a commutative kei we get the result immediately, 
since our assumption implies $x*z=y*z$ and that forces equality $x=y$.
Let us consider the case $k>3$.
Since the relation $r_k$ holds for all elements, we have:
$$z=x*z*x*\ldots *z*x\ \textrm{(with $k$ letters on the right)} ,$$
which is equivalent to
$$z*x*z=x*z*x*\ldots *z*x*z*x\ \textrm{(with $k-2$ letters on the right)}$$
and $$z=y*z*y*\ldots *z*y,$$
equivalent to
$$z*y*z=y*z*y*\ldots *z*y*z*y.$$
We can replace the initial assumption with:
$$ z*x*z=z*y*z $$
and use relation $r_k$ to change this equation to:
$$ x*z*x*\ldots *z*x*z*x=y*z*y*\ldots *z*y*z*y;$$
$$ x*z*x*\ldots *z*(z*x)=y*z*y*\ldots *z*y*z*y;$$
$$ x*z*x*\ldots *z*(z*y)=y*z*y*\ldots *z*y*z*y;$$
$$ x*z*x*\ldots *z*y*z*y=y*z*y*\ldots *z*y*z*y.$$
Now we can cancel last four letters from both sides of the equation.
This reduction is repeated until we arrive to $x=y$ or $x*z*x=y*z*y$. In the latter case we use the first quandle axiom to write:
$$ x*x*z*x=y*y*z*y; $$
$$ x*(z*x)=y*(z*y)$$
and reduce one more time to obtain $x=y$.
\end{proof}
For every quandle $Q$ we can consider its operator group, $\it{Op} (Q)$, generated by automorphisms $f_x \colon Q \to Q$, defined by $f_x(y)=y*x$.
Behaviorally equivalent elements of $Q$ define equal elements in $\it{Op} (Q)$.
On the other hand, if there are no behaviorally equivalent elements in $Q$ then the map
$x \mapsto f_x$ is an injection and $Q$ is isomorphic to a union of conjugacy classes (of images of generators of $Q$)
in $\it{Op} (Q)$ (see also \cite{Joy}).
\begin{corollary} \label {emb}
Any kei $Q$ satisfying the universal relation $r_k$, for some odd $k$, 
embeds into the conjugation quandle of its operator group, 
$\textit{Conj}(\textit{Op} (Q))$ (with quandle operation $f*g=g^{-1}fg$).
\end{corollary}
Absence of behaviorally equivalent elements enables us to prove the following theorem.
\begin{theorem} \label{3^k}
The order of a finite commutative kei $Q$ is a power of $3$.
\end{theorem}
\begin{proof}
Let $a$, $b$ denote any two elements of $Q$ and let $P$ be the three element
subquandle $\{a,\,b,\,a*b\}$. From the right distributivity property it 
follows that for any element $x\in Q$, the set $S:=P*x=\{a*x,\,b*x,\,(a*b)*x\}$
is also a subquandle. We ask the following question: 
what other elements of $Q$ send $P$ to $S$ ?
There can be at most three such operators, 
sending $a$ to $a*x$, $b*x$ or $(a*b)*x$ (here we use the fact that 
two operators $x$, $y\in Q$ acting in the same way on one element 
are the same). Using lemma \ref{strong} we can find them easily:
\begin{enumerate}
\item[(1)] $a*(a*x)=\mathbf{x}$
\item[(2)] $a*(b*x)=b*x*a=\mathbf{x*b*a}$;
\item[(3)] $a*(a*b*x)=a*b*x*a=b*a*x*a=b*(x*a)=\mathbf{x*a*b}$.
\end{enumerate}
From the left distributivity it follows that the set of these three operators 
is a subquandle.
We still need to check that operators (2) and (3) send $b$ and $a*b$ to $S$:
$$b*(x*b*a)=x*b*a*b=x*(a*b)=(a*b)*x;$$
$$(a*b)*(x*b*a)=b*a*(x*b*a)=b*a*a*b*x*b*a=b*x*b*a=x*b*b*a=a*x;$$
$$b*(x*a*b)=x*a*b*b=x*a=a*x;$$
$$(a*b)*(x*a*b)=a*b*b*a*x*a*b=a*x*a*b=x*a*a*b=x*b=b*x.$$
In this way we obtain a partition of $Q$ into 3-element disjoint 
subquandles of the form $\{x,\,x*b*a,\,x*a*b\}$, in which two elements 
belong to the same triple if they send $P$ to the same subquandle. 
This relation between elements is an equivalence relation but not 
a congruence ($u\sim v$, $s\sim t$ does not imply $(u*s)\sim (v*t)$), 
so we cannot simply form a quotient quandle. Instead, we define 
a natural quandle operation, $\hat *$, on triples:
$$(x,\,x*b*a,\,x*a*b)\,\hat{*}\,(y,\,y*b*a,\,y*a*b)=(x*y,\,x*y*b*a,\,x*y*a*b).$$
The set of such triples, with operation $\hat{*}$, forms a commutative 
kei that is three times smaller than the original kei, $Q$. 
Thus we can use the inductive argument to conclude that the size of $Q$ 
is a power of $3$.
\end{proof} 
\subsection{$\mathbf{\bar Q(4,3)}$ has $81$ elements}
It was shown by M. Takasaki \cite{Tak} that $\bar Q(2,n)$ is isomorphic 
to the dihedral quandle $Z_n$ and $\bar Q(3,3)$ is isomorphic 
to $Z_3 \oplus Z_3$.
Here we give a description of $\bar Q(4,3)$.\\
T.Ohtsuki wrote a computer program which helps to analyze the commutative kei. 
Using this program he found that $\bar Q(4,3)$ has $81$ elements. 
A different computation, involving operator group of the quandle, 
was made by the first author.
Here we follow, in a crucial point, Ohtsuki's approach to obtain 
a computer free proof.  
\begin{theorem} \label{main}
$\bar Q(4,3)$ has 81 elements.
\end{theorem}
As noted by Takasaki, every element of the kei can be written 
in a left-normed form (usually not uniquely).
For example, in 
$\bar Q(4,3)$, $(a*b)*(c*d)= a*b*c*d*c=a*b*d*c*d=b*a*c*d*c= 
b*a*d*c*d=c*d*a*b*a=c*d*b*a*b=d*c*a*b*a=d*c*b*a*b$. The length of the 
kei element $w$ is the length of the shortest left-normed word, 
in the generators of kei, representing $w$.
\begin{lemma} \begin{enumerate} \label{l1}
\item[(i)] Every element of $\bar Q(4,3)$, in a generating set 
$\{a,b,c,d\}$, is of length at most 7.
\item [(ii)] There are (at most) $8$ elements in $\bar Q(4,3)$ of length $7$ and 
they have representatives:\\
$a*b*c*d*b*c*d$, $a*b*d*c*b*d*c$,\\
$b*a*c*d*a*c*d$, $b*a*d*c*a*d*c$,\\
$c*a*b*d*a*b*d$, $c*a*d*b*a*d*b$,\\
$d*a*b*c*a*b*c$, $d*a*c*b*a*c*b$,
\end{enumerate}
\end{lemma}
\begin{proof} 
We use brackets [\ ] to stress  
which group of letters our properties are used on.
The bracket [\ ] (unlike $(\ )$) does not change the left-normed convention.
Let $\{x_0,x_1,x_2,x_3\}=\{a,b,c,d\}$. We have the following 
identities in $\bar Q(4,3)$.
\begin{enumerate}
\item[(1)] $x_0*x_1*x_2*x_0=x_0*x_2*x_1$.\\
This is the case because 
$x_0*x_1*x_2*x_0= x_1*x_0*x_2*x_0=x_1*(x_0*x_2)=x_0*x_2*x_1$.
\item[(2)] $w*x_0*x_1*x_0=w*x_1*x_0*x_1$.\\
It is the case because $w*x_0*x_1*x_0=w*(x_1*x_0)=w*(x_0*x_1)=w*x_1*x_0*x_1$.
\item[(3)] $x_0*x_1*x_2*x_3*x_2= (x_0*x_1)*(x_2*x_3)= 
(x_2*x_3)*(x_0*x_1)=x_2*x_3*x_0*x_1*x_0$.
\item[(4)] $x_0*x_1*x_2*x_3*x_2*x$ is 
reducible to a word of length $4$ for $x=x_i$ 
$i=0,1,2,3$, for example $x_0*x_1*x_2*x_3*x_2*x_1=(x_0*x_1)*(x_2*x_3)*x_1=
x_2*x_3*x_1*x_0$.
\item[(5)] $x_0*x_1*x_2*x_3*x_0*x_1 = x_0*x_1*x_3*x_2*x_0$.\\
It is the case because $x_0*x_1*x_2*x_3*x_0*x_1*x_0 = 
((x_0*x_1)*(x_2*x_3))*x_2*(x_0*x_1) = (x_2*x_3)*(x_0*x_1*x_2)= (x_0*x_1*x_2)*(x_2*x_3)=
x_0*x_1*x_3*x_2$ as required. 
\item[(6)] $x_0*x_1*x_2*x_3*x_0*x_3 = x_0*x_2*x_1*x_3*x_0$.\\
It is the case because  $x_0*x_1*x_2*x_3*x_0*x_3= x_1*x_0*x_2*x_3*x_0*x_3=^{(1)}
x_1*x_0*x_2*x_0*x_3*x_0=^{(2)} x_0*x_2*x_1*x_3*x_0$.
\item[(7)] $x_0*x_1*x_2*x_3*x_1*x_2= x_3*x_2*x_1*x_0*x_2*x_1$.\\
It is the case because $x_0*x_1*x_2*x_3*x_1*x_2= ((x_0*x_1)*(x_2*x_3))*(x_1*x_2)=
((x_3*x_2)*(x_1*x_0))*(x_2*x_1)= x_3*x_2*x_1*x_0*x_2*x_1$.
\item[(8)] $x_0*x_1*x_2*x_3*x_1*x_2*x_3=x_0*x_2*x_3*x_1*x_2*x_3*x_1=x_0*x_3*x_1*x_2*x_3*x_1*x_2$. 
The equalities hold because\\
$x_0*x_1*x_2*x_3*x_1*[x_2*x_3*x_2]=^{(2)} (x_0*x_1)*x_2*x_3*x_1*x_3*x_2*x_3=
(x_1*x_0*x_2*x_3*x_1*x_3)*x_2*x_3=^{(6)}$ \\
$(x_1*x_2*x_0*x_3*x_1*x_2)*x_3 =^{(5)} x_1*x_2*x_3*x_0*x_1*x_3=$\\
 $ x_2*x_1*x_3*x_0*x_1*x_3 =^{(7)} x_0*x_3*x_1*x_2*x_3*x_1$ as required.
\item[(9)] $x_0*x_1*x_2*x_3*x_1*x_2*x_3*x_0 = x_0*x_1*x_3*x_2*x_1*x_3*x_2$. This 
equality is the most difficult and allows us to complete Lemma \ref{l1}.
In the proof we follow Ohtsuki's analysis of his 
computer computation. He noticed that the key computation 
is to use the commutative identity:
$(x_0*x_1*x_2*x_3*x_1*x_2*x_3*x_0)*x_2=x_2*(x_0*x_1*x_2*x_3*x_1*x_2*x_3*x_0)$
and to show that the last expression can be reduced to $x_0*x_1*x_3*x_2*x_1*x_3$ 
by properties (1)-(8).
\end{enumerate}
\end{proof}
\begin{lemma} \label{l2}
$x_2*(x_0*x_1*x_2*x_3*x_1*x_2*x_3*x_0)=x_0*x_1*x_3*x_2*x_1*x_3$
\end{lemma}
\begin{proof} For improved clarity, we omit * in the presentation of words in this proof. Using the identity $wxyx=w(yx)=w(xy)=wyxy$ for seven times we obtain
$x_2(x_0x_1x_2x_3x_1x_2x_3x_0)= 
x_2x_0x_3x_2x_1x_3x_2x_1x_0x_1x_2x_3x_1x_2x_3x_0$. Next we use identities 
(1)-(8) for several times to get\\
$(x_2x_0x_3x_2)x_1x_3x_2x_1x_0x_1x_2x_3x_1x_2x_3x_0=^{(1)}$\\  
$x_2x_3x_0x_1x_3x_2x_1x_0x_1x_2x_3x_1x_2x_3x_0
= (x_3x_2x_0x_1x_3x_2)x_1x_0x_1x_2x_3x_1x_2x_3x_0 =^{(5)} $\\ 
$x_3x_2x_1x_0x_3[x_1x_0x_1]x_2x_3x_1x_2x_3x_0=^{(2)} 
(x_3x_2x_1x_0x_3x_0)x_1x_0x_2x_3x_1x_2x_3x_0 =^{(6)} $\\
$(x_3x_1x_2x_0x_3 x_1)x_0x_2x_3x_1x_2x_3x_0 =^{(5)}$ \\ 
$(x_3x_1)x_0x_2x_3x_0x_2x_3x_1x_2x_3x_0= 
(x_1x_3x_0x_2x_3x_0)x_2x_3x_1x_2x_3x_0 =^{(7)}$\\
$ x_2x_0x_3x_1x_0[x_3x_2x_3]x_1x_2x_3x_0=^{(2)}
 (x_2x_0)x_3x_1x_0x_2x_3x_2x_1x_2x_3x_0= $\\ 
$(x_0x_2x_3x_1x_0x_2)x_3x_2x_1x_2x_3x_0 =^{(5)}  
(x_2x_0)x_1x_3x_0x_3x_2x_1x_2x_3x_0 = $\\ 
$(x_0x_2x_1x_3x_0x_3)x_2x_1x_2x_3x_0 =^{(6)} 
x_0x_1x_2x_3x_0[x_2x_1x_2]x_3x_0 =^{(2)}$ \\
$ (x_0x_1x_2x_3x_0x_1)x_2x_1x_3x_0 =^{(5)} (x_0x_1x_3x_2x_0x_2)x_1x_3x_0 
=^{(6)} $\\ 
$ (x_0x_3)x_1x_2x_0x_1x_3x_0 = (x_3x_0x_1x_2x_0x_1)x_3x_0 =^{(7)}$\\ 
$x_2x_1x_0x_3x_1[x_0x_3x_0] =^{(2)} (x_2x_1)x_0x_3x_1x_3x_0x_3=$\\
$(x_1x_2x_0x_3x_1x_3)x_0x_3 =^{(6)} (x_1x_0x_2x_3x_1x_0)x_3=^{(5)} 
(x_1x_0)x_3x_2x_1x_3 = x_0x_1x_3x_2x_1x_3$ as required.
\end{proof} 
We proved, in Lemma \ref{l1}, that $\bar Q(4,3)$ is finite, 
but in fact we can easily build, using Lemmas \ref{l1}, \ref{l2} 
and their proofs, the multiplication table of $\bar Q(4,3)$, 
with 81 elements. We should still argue that $\bar Q(4,3)$ is not smaller. 
One argument, very laborious and good for the computer, 
is that we can use all relations of the commutative kei and no reduction 
will be found. More sophisticated argument uses the kei epimorphism 
$p: \bar Q(4,3) \to Z_3^3$.
$p$ is given on generators of $\bar Q(4,3)$ by:\\
$p(a)=(0,0,0)$, $p(b)=(1,0,0)$, $p(c)=(0,1,0)$ and $p(d)=(0,0,1)$.\\ 
For example, we have $p(a*b)=(2,0,0)$, $p(a*c)=(0,2,0)$, $p(a*d)=(0,0,2)$,
$p(b*c)=(2,2,0)$, $p(b*d)= (2,0,2)$, $p(c*d)=(0,2,2)$.\\
From theorem \ref{3^k} it follows that it is enough to prove 
that $p$ is not a monomorphism.
We notice that $p((a*b)*(c*d))=(1,1,1)=p((a*c)*(b*d))=p((a*d)*(b*d))$.\\
However, in $\bar Q(4,3)$ we have the inequality 
$(a*b)*(c*d)\neq (a*c)*(b*d)$.\\
To prove this inequality we use another homomorphism  
$q: \bar Q(4,3) \to B(4,3)$  being the identity on generators.
We have to check whether\\$q((a*b)*(c*d))=q((a*c)*(b*d))$.\\
We have $q((a*b)*(c*d))= q(c*d)(q(a*b))^{-1}q(c*d)=
dc^{-1}d(b^{-1}ab^{-1})dc^{-1}d,$ similarly 
$q((a*c)*(b*d))= db^{-1}d(c^{-1}ac^{-1})db^{-1}d$.\\Thus we have to check 
that\\$q((a*b)*(c*d))(q((a*c)*(b*d)))^{-1}= 
dc^{-1}d(b^{-1}ab^{-1})dc^{-1}bd^{-1}ca^{-1}cd^{-1}bd^{-1}$\\
is not equal to $1$ in $B(4,3)$. 
After conjugating it by $d^{-1}$ we 
reduce it to the question
$$c^{-1}db^{-1}ab^{-1}dc^{-1}bd^{-1}ca^{-1}cd^{-1}b \neq 1\ ?$$ 
This inequality was confirmed by GAP but 
Mietek D{\c a}bkowski also checked it by hand, using the lower central 
series of the Burnside group.
It follows that $p$ is not a monomorphism and $\bar Q(4,3)$ has 
exactly $81$ elements.\\
\textbf{Remark (alternative approach)}. 
For every quandle $Q$ we can consider its associated group, $\textit{As}(Q)$, defined as the quotient $F(Q)/K$, where $F(Q)$ denotes the free group on elements of $Q$ and $K$ is the normal subgroup generated by the words $(x*y)y^{-1}x^{-1}y$, where $x$, $y\in Q$.
The operator group, $\textit{Op}(Q)$, is the quotient of the associated 
group of $Q$ (see for example \cite{F-R}). From this fact and the theorem 
of Winker (which we restate slightly modified to match our notation) 
it follows that the operator group of $\bar Q(4,3)$ is generated by 
the images of generators $a$, $b$, $c$, $d$.
\begin{theorem}[\cite{Win}, theorem 5.1.7]
Let $\{ S\,\mid\, R\}$ be the presentation of the quandle $Q$.
Then the group $\textit{As}(Q)$ has a presentation $\{\bar{S}\,\mid\,\bar{R}\}$, where $\bar{S}=\{\bar{s}\,\colon\, s\in S\}$ and
$\bar{R}=\{\bar{r}=\bar{s}\,\colon\, r=s$ is a relation in $R\}$.
Here $\bar{r}$ denotes the group element obtained by replacing $u*v$ by 
$\bar{v}^{-1}\bar{u}\bar{v}$.
\end{theorem}
As we noted previously (see Corollary \ref {emb} and the comment preceding it), 
$\bar Q(4,3)$ embeds into its operator group as conjugates of generators. 
For simplicity we will use the same notation for the quandle elements 
and their images in the operator group. From the second kei axiom 
 it follows that the squares of generators (and therefore also the squares 
of conjugates of generators) are equal to the identity in 
$\textit{Op}(\bar Q(4,3))$. 
Our commutative relation, $x*y=y*x$, which is true for all 
elements of $\bar Q(4,3)$,
becomes $yxy=xyx$ (or $xyxyxy=1$) in $\textit{Op}(\bar Q(4,3))$, 
where $x$ and $y$ belong to conjugacy classes of generators. 
Now it follows that $\bar Q(4,3)$ embeds into a (possibly bigger) 
group with presentation:
$$\{ a,\, b,\, c,\, d\,\mid\, a^2=b^2=c^2=d^2=1,\, xyxyxy=1 \} ,$$\ 
where $x$ and $y$ are any conjugates of the generators $a$, $b$, $c$, $d$.
We computed (using GAP) that the order of this group 
is $118098=2\cdot 3^{10}$, and the number of elements 
in conjugacy classes of $a$, $b$, $c$, $d$ is $81$, which is also the size of 
$\bar Q(4,3)$ (elements of these conjugacy classes 
form a $4$-generator commutative kei with conjugation as a quandle 
operation, therefore their number cannot exceed the order 
of the free kei $\bar Q(4,3)$). 
\subsection{$\mathbf{\bar Q(4,3)}$ as an extension of 
$\mathbf{Z_3^3}$ by $\mathbf{Z_3}$}
We show that $\bar Q(4,3)$ can be represented as a quandle 
$(Z_3 \times Z_3^3,\hat{*})$,
where the operation $\hat{*}$ is defined on the set $Z_3 \times Z_3^3$ by: 
$$(a_1,x_1)\hat{*} (a_2,x_2)=(a_1*a_2 + c(x_1,x_2),x_1*x_2).$$
Here $+$ is addition in $Z_3$ and $c(x_1,x_2)\colon Z_3^3\times Z_3^3 
\to Z_3$ is a function that must satisfy the following conditions 
coming from the kei axioms.
\begin{enumerate}
\item[(i)] The axiom $(a,x)\hat{*} (a,x)= (a,x)$ leads to 
$c(x,x)=0$.
\item[(ii)] The condition $((a,x)\hat{*}(b,y))\hat{*}(b,y)=
(a,x)$ leads to $c(x*y,y)=c(x,y)$.
\item[(iii)] The distributivity property
$$((a_1,x_1)\hat{*}(a_2,x_2))\hat{*}(a_3,x_3) = 
((a_1,x_1)\hat{*}(a_3,x_3))\hat{*}((a_2,x_2)\hat{*}(a_3,x_3))$$ leads to the following, 
after first computing the left and the right side of the above equation.\\
$L=(a_1*a_2 + c(x_1,x_2),x_1*x_2)\hat{*}(a_3,x_3) =
(a_1*a_2*a_3 - c(x_1,x_2) + c(x_1*x_2,x_3),x_1*x_2*x_3)$.\\
$R= (a_1*a_3 + c(x_1,x_3),x_1*x_3)\hat{*}
    (a_2*a_3 + c(x_2,x_3),x_2*x_3)=
((a_1*a_3)*(a_2*a_3)+ 2c(x_2,x_3)-c(x_1,x_3) + 
c(x_1*x_3,x_2*x_3),(x_1*x_3)*(x_2*x_3))$\\
From this we get:
$c(x_1*x_3,x_2*x_3) - c(x_1*x_2,x_3) = - c(x_1,x_2) -2c(x_2,x_3) +c(x_1,x_3)$.
Taking into account that we work modulo $3$, we get:
$$c(x_1*x_3,x_2*x_3) - c(x_1*x_2,x_3) = - c(x_1,x_2) + c(x_2,x_3) +c(x_1,x_3)$$
\item[(iv)] The condition that our kei is commutative leads to $c(x,y)=c(y,x)$.
\end{enumerate}
The condition (iii) is what makes function $c$ to be a twisted $2$-cocycle in the second quandle cohomology group of $Z_3^3$ with coefficients in $Z_3$. Twisted quandle (co)homology theory was introduced in \cite{C-E-S}. The authors described there a general method of obtaining a new quandle from a given quandle $X$ and Alexander quandle $A$, using a twisted $2$-cocycle $\phi$.\\
Such constructions, including the one we are describing, are called Alexander extensions of $X$ by $(A,\phi)$.\\
An example of the function $c$ satisfying all of the above conditions, 
is presented below. We need to order (assign numbers to) the elements of 
$Z_3^3$ in order to describe the matrix defining cocycle $c$:\\
1.(0,0,0); 2.(0,0,1); 3.(0,0,2); 4.(0,1,0); 5.(0,1,1); 6.(0,1,2); 7.(0,2,0);\\
8.(0,2,1); 9.(0,2,2); 10.(1,0,0); 11.(1,0,1); 12.(1,0,2); 13.(1,1,0); 
14.(1,1,1);\\ 15.(1,1,2); 16.(1,2,0); 17.(1,2,1); 18.(1,2,2); 19.(2,0,0); 
20.(2,0,1); 21.(2,0,2);\\ 22.(2,1,0); 23.(2,1,1); 24.(2,1,2); 
25.(2,2,0); 26.(2,2,1); 27.(2,2,2).\\
In the matrix $M$, the entry $m_{ij}$ is equal to $c(i,j)$, the value 
of the cocycle $c$ on elements numbered with $i$, $j$.
{\small
\begin{displaymath}
M=\left[ 
\begin{tabular}{l @{\hspace{0.2 cm}} l @{\hspace{0.2 cm}} l @{\hspace{0.2 cm}} l @{\hspace{0.2 cm}}
l @{\hspace{0.21 cm}} l @{\hspace{0.21 cm}} l @{\hspace{0.21 cm}} l @{\hspace{0.21 cm}}
l @{\hspace{0.21 cm}} l @{\hspace{0.21 cm}} l @{\hspace{0.21 cm}} l @{\hspace{0.21 cm}}
l @{\hspace{0.21 cm}} l @{\hspace{0.21 cm}} l @{\hspace{0.21 cm}} l @{\hspace{0.21 cm}}
l @{\hspace{0.21 cm}} l @{\hspace{0.21 cm}} l @{\hspace{0.21 cm}} l @{\hspace{0.21 cm}}
l @{\hspace{0.21 cm}} l @{\hspace{0.21 cm}} l @{\hspace{0.21 cm}} l @{\hspace{0.21 cm}}
l @{\hspace{0.21 cm}} l @{\hspace{0.21 cm}} l}
0&1&1&0&1&2&0&2&1&0&1&2&0&1&2&0&1&2&0&2&1&0&2&1&0&2&1\\
1&0&1&0&1&2&2&1&0&2&0&1&1&2&0&0&1&2&1&0&2&2&1&0&0&2&1\\ 
1&1&0&1&2&0&2&1&0&2&0&1&0&1&2&1&2&0&0&2&1&2&1&0&1&0&2\\ 
0&0&1&0&0&0&0&1&0&2&0&1&1&2&0&0&1&2&0&2&1&1&0&2&2&1&0\\ 
1&1&2&0&0&0&2&1&1&0&1&2&1&2&0&2&0&1&1&0&2&0&2&1&2&1&0\\ 
2&2&0&0&0&0&2&2&0&1&2&0&1&2&0&1&2&0&2&1&0&2&1&0&2&1&0\\ 
0&2&2&0&2&2&0&0&0&0&1&2&1&2&0&2&0&1&1&0&2&0&2&1&2&1&0\\ 
2&1&1&1&1&2&0&0&0&1&2&0&1&2&0&1&2&0&0&2&1&0&2&1&0&2&1\\ 
1&0&0&0&1&0&0&0&0&2&0&1&1&2&0&0&1&2&2&1&0&0&2&1&1&0&2\\ 
0&2&2&2&0&1&0&1&2&0&0&0&1&2&0&1&0&2&0&2&2&0&2&1&2&1&0\\ 
1&0&0&0&1&2&1&2&0&0&0&0&0&1&2&2&1&0&0&0&1&0&2&1&2&1&0\\
2&1&1&1&2&0&2&0&1&0&0&0&2&0&1&0&2&1&1&2&1&0&2&1&2&1&0\\ 
0&1&0&1&1&1&1&1&1&1&0&2&0&0&0&1&2&0&1&1&1&1&1&1&0&0&1\\ 
1&2&1&2&2&2&2&2&2&2&1&0&0&0&0&0&1&2&2&2&2&2&2&2&1&2&1\\ 
2&0&2&0&0&0&0&0&0&0&2&1&0&0&0&2&0&1&0&0&0&0&0&0&0&2&2\\ 
0&0&1&0&2&1&2&1&0&1&2&0&1&0&2&0&0&0&0&1&2&0&1&0&2&0&1\\ 
1&1&2&1&0&2&0&2&1&0&1&2&2&1&0&0&0&0&2&0&1&2&1&1&1&2&0\\ 
2&2&0&2&1&0&1&0&2&2&0&1&0&2&1&0&0&0&1&2&0&2&2&0&0&1&2\\ 
0&1&0&0&1&2&1&0&2&0&0&1&1&2&0&0&2&1&0&0&0&2&1&0&2&0&1\\ 
2&0&2&2&0&1&0&2&1&2&0&2&1&2&0&1&0&2&0&0&0&0&2&1&1&2&0\\ 
1&2&1&1&2&0&2&1&0&2&1&1&1&2&0&2&1&0&0&0&0&1&0&2&0&1&2\\
0&2&2&1&0&2&0&0&0&0&0&0&1&2&0&0&2&2&2&0&1&0&0&0&2&1&0\\ 
2&1&1&0&2&1&2&2&2&2&2&2&1&2&0&1&1&2&1&2&0&0&0&0&0&2&1\\ 
1&0&0&2&1&0&1&1&1&1&1&1&1&2&0&0&1&0&0&1&2&0&0&0&1&0&2\\ 
0&0&1&2&2&2&2&0&1&2&2&2&0&1&0&2&1&0&2&1&0&2&0&1&0&0&0\\ 
2&2&0&1&1&1&1&2&0&1&1&1&0&2&2&0&2&1&0&2&1&1&2&0&0&0&0\\ 
1&1&2&0&0&0&0&1&2&0&0&0&1&1&2&1&0&2&1&0&2&0&1&2&0&0&0  
\end{tabular}
\right] 
\end{displaymath} }
The isomorphism between $\bar Q(4,3)$ and $(Z_3 \times Z_3^3,\hat{*})$ follows from the facts:
\begin{enumerate}
\item[(i)] $(Z_3 \times Z_3^3,\hat{*})$ satisfies kei axioms and relation $r_3$;
\item[(ii)] $(Z_3 \times Z_3^3,\hat{*})$ has 81 elements;
\item[(iii)] $(Z_3 \times Z_3^3,\hat{*})$ is generated by four elements:\\ 
(1,0,0,0), (0,1,0,0), (0,0,1,0), (0,0,0,1) (for example it cannot be isomorphic to 
$\textit{Core}(Z_3^4)$, which has five generators as a kei).
\end{enumerate}
\section{{$\mathbf{\bar Q(3,4)}$} has $96$ elements.}
\begin{figure}
\begin{center}
\includegraphics[height=14 cm]{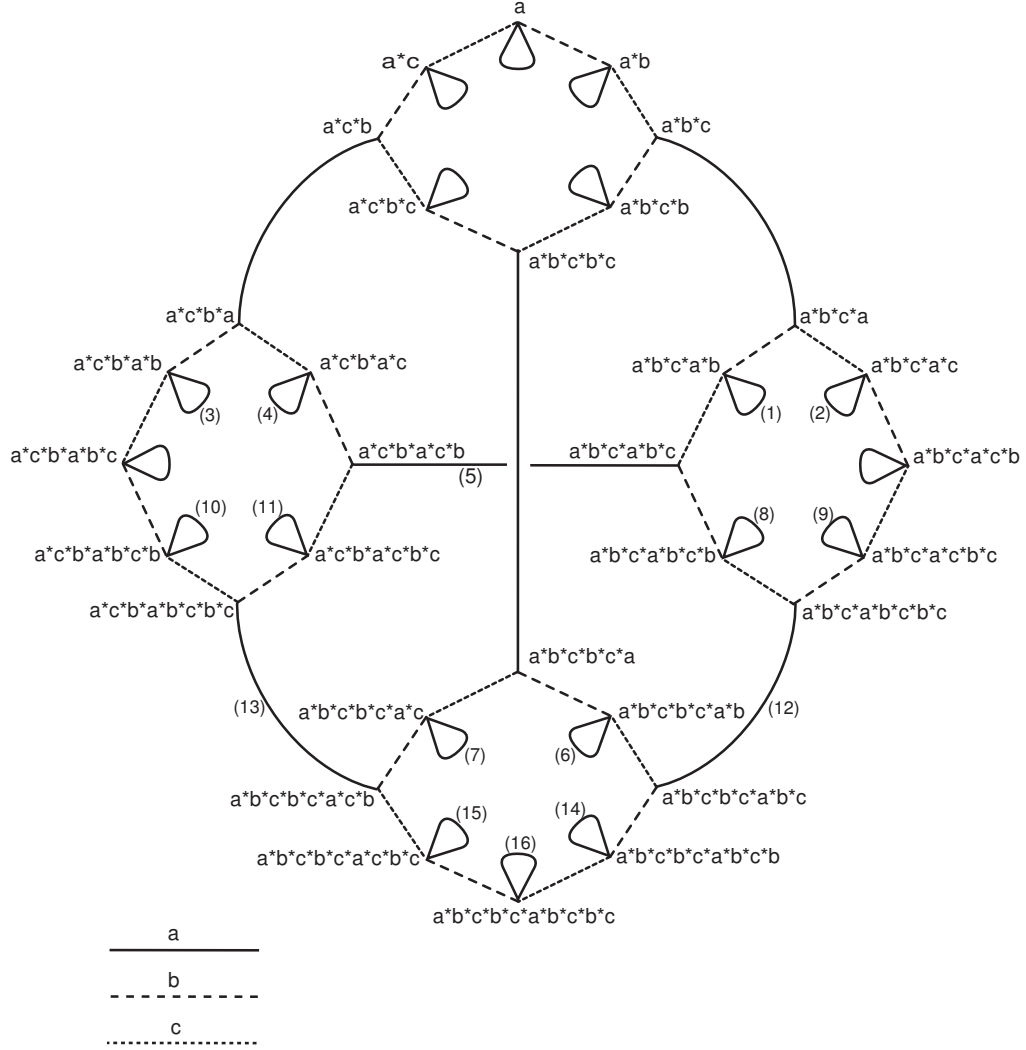}
\caption{One of the components in the Cayley diagram of 
$\bar Q(3,4)$.\label{q34}}
\end{center}
\end{figure}

The primary examples of keis satisfying the fourth Burnside 
relation, $x=x*y*x*y$, are the dihedral kei $Z_4$, its direct sums, 
and the fourth Burnside groups and their quotients (with the core 
operation $x*y=yx^{-1}y)$.\\
To get the lower bound on the order of $\bar Q(3,4)$, 
we can consider the group with presentation:
$$\{ a,\,b,\,c\,\mid\, a^2=b^2=c^2=1,\,(xy)^4=1\} ,$$
where $x$ and $y$ are any conjugates of generators $a$, $b$, $c$.
We checked (with the help of GAP) that this group has $8192=2^{13}$ 
elements and the size of the union of conjugacy classes of generators 
is $96$. The elements of these conjugacy classes form a $3$-generator 
quandle (with conjugation as operation $*$) satisfying 
the relation $x=x*y*x*y$.
Thus the order of the free kei $\bar Q(3,4)$ cannot be less than $96$.
This time we cannot use the Lemma \ref{belements} to obtain 
the upper bound for the size of $\bar Q(3,4)$ (compare with the 
remark after the proof of Theorem \ref{main}), so instead we will 
build a Cayley diagram for this quandle.
This diagram has $96$ vertices, therefore $\bar Q(3,4)$ has order $96$.
Below we calculate some relations needed to build
such a diagram. Again, we use brackets $[\;]$ 
to stress which parts of words properties
of $\bar Q(3,4)$ are used on. None of these properties can replace the first 
letter in the left normed representatives of words in $\bar Q(3,4)$. 
For example a left normed word starting with $a$ never equals to the word 
starting with $b$. It follows that the diagram will consist of three 
disjoint parts that look the same when viewed 
as graphs (see Figure \ref{q34}).
Here we prove the most difficult relations in the Cayley graph, their numbers correspond to the numbers included in the Figure \ref{q34}.
The vertices of the Cayley graph represent the elements of $\bar Q(3,4)$. The solid arcs represent multiplication from the right by the generator $a$; two kinds of dashed arcs denote multiplication by respectively $b$ and $c$.
\begin{enumerate}
\item[(1)] We need to prove the relation $a*b*c*a*b*a=a*b*c*a*b$ (corresponding to a loop at a vertex representing element $a*b*c*a*b$):
$[a*b]*c*a*b*a=a*[b*a*c*a*b]*a=a*(c*a*b)*a=a*(c*a*b)=[a*b*a]*c*a*b=a*b*c*a*b$ as wanted;
\item[(2)] $[a*b]*c*a*c*a=a*b*[a*c*a*c*a]=a*b*[(a*c*a)]=$\\$a*b*[(a*c)]=a*b*c*a*c$;
\item[(3)] Similar to (2);
\item[(4)] Similar to (1);
\item[(5)] We have to check that $a*b*c*a*b*c*a=a*c*b*a*c*b$ or equivalently that
$a=a*c*b*a*c*b*a*c*b*a*c*b$:\\
$[a*c]*b*a*c*b*a*c*b*a*c*b=a*[c*a*b*a*c]*b*a*c*b*a*c*b=[a*(b*a*c)]*b*a*c*b*a*c*b=$
$a*[(b*a*c)]*a*b*a*c*b*a*c*b=a*[c*a*b*a*c*a*b*a*c]*b*a*c*b=$
$a*[(c*(b*a)*c)]*b*a*c*b=a*[(c*(b*a))]*b*a*c*b=a*a*b*a*c*a*[*b*a*b*a]*c*b=$
$[a*a]*b*a*c*[a*a]*b*a*b*c*b=[a*b*a]*c*b*a*b*c*b=a*[b*c*b*a*b*c*b]=a*(a*(c*b))=a$;
\item[(6)] We prove that $a*b*c*b*c*a*b*a=a*b*c*b*c*a*b$:
$a*b*c*b*c*a*b*a=a*b*c*b*c*b*[b*a*b*a]=a*[b*c*b*c*b]*a*b*a*b=
[a*c*b*c*a]*b*a*b=a*[c*b*c*b]*a*b=a*b*c*b*c*a*b$;
\item[(7)] $a*[b*c*b*c]*a*c*a=a*c*b*c*b*a*c*a$ and then like in (6) (the roles of $b$ and $c$ are exchanged);
\item[(8)] We need $a*b*c*a*b*c*b=a*b*c*a*b*c*b*a$ or $a*b*c*a*b*c*b*a*b*c*b*a*c*b=a$,\\
$[a*b]*c*a*b*c*b*a*b*c*b*a*c*b=a*[b*a*c*a*b]*c*b*a*b*c*b*a*c*b=
[a*(c*a*b)]*c*b*a*b*c*b*a*c*b=a*[(c*a*b)]*a*c*b*a*b*c*b*a*c*b=
a*b*a*c*a*b*a*c*b*a*b*c*b*a*c*b=
a*b*a*c*a*b*a*c*a*[a*b*a*b]*c*b*a*c*b=a*[b*a*c*a*b*a*c*a*b]*a*b*a*c*b*a*c*b=
a*[(b*(c*a)*b)]*a*b*a*c*b*a*c*b=
a*[(b*a*c*a)]*a*b*a*c*b*a*c*b=[a*a]*c*a*b*a*c*[a*a]*b*a*c*b*a*c*b=
[a*c*a]*b*a*c*b*a*c*b*a*c*b=a*c*b*a*c*b*a*c*b*a*c*b=^{(5)}a$;
\item[(9)] $a*b*c*a*c*b*c=a*b*c*a*c*b*c*a$ is equivalent to \\
$a*b*c*a*c*b*c*a*c*b*c*a*c*b=a$, \\
$[a*b*c*a*c*b]*c*a*c*b*c*a*c*b=a*b*c*a*c*b*[a*c*a*c]*b*c*a*c*b=
a*b*c*a*c*b*c*a*c*a*b*c*a*c*b=a*[b*c*a*c*b*c*a*c*b]*b*a*b*c*a*c*b=
a*[(b*(a*c)*b)]*b*a*b*c*a*c*b=[a*c*a]*c*b*c*a*c*b*a*b*c*a*c*b=
a*[c*c]*b*c*a*c*b*a*b*c*a*c*b=a*[b*c*a*c*b*a*b*c*a*c*b]=
a*(a*b*c*a*c*b)=a*(a*(a*c*b))=a$;
\item[(10)] Like in (9) with $b$ and $c$ interchanged;
\item[(11)] Like in (8) with $b$ and $c$ interchanged;
\item[(12)] $a*b*c*a*b*c*b*c*a=a*b*c*b*c*a*b*c$ is equivalent to\\
$a*b*c*b*c*a*b*c*a*c*b*c*b*a*c*b=a$, \\
$[a*b]*c*b*c*a*b*c*a*c*b*c*b*a*c*b=a*[b*a*c*b*c*a*b]*c*a*c*b*c*b*a*c*b=
[a*(b*c*a*b)]*c*a*c*b*c*b*a*c*b=a*[(b*c*a*b)]*a*c*a*c*b*c*b*a*c*b=
a*b*a*c*b*c*a*b*a*c*a*[c*b*c*b]*a*c*b=
[a*b*a]*c*b*c*a*b*a*c*a*b*c*b*c*a*c*b=[a*b*c*b]*c*a*b*a*c*a*b*c*b*c*a*c*b=
a*b*c*[b*a*c*a*b*a*c*a*b]*c*b*c*a*c*b=a*b*c*[(b*(c*a)*b)]*c*b*c*a*c*b
a*b*c*a*c*a*b*a*c*a*c*b*c*a*c*b=a*b*c*a*c*a*b*a*b*[b*c*a*c*b*c*a*c*b]=
[a*b]*c*a*c*a*b*a*b*c*a*c*b*c*a*c=a*b*[a*c*a*c*a]*b*a*b*c*a*c*b*c*a*c=
a*[b*c*a*c*b*a*b*c*a*c*b]*c*a*c=a*(a*(a*c*b))c*a*c=a*c*a*c=a$;
\item[(13)] Follows from (12);
\item[(14)] Instead of $a*b*c*b*c*a*b*c*b*a=a*b*c*b*c*a*b*c*b$ we consider
$a=a*b*c*b*c*a*b*c*b*a*b*c*b*a*c*b*c*b$,\\
$a*b*c*b*c*[a*b*c*b*a*b*c*b*a]*c*b*c*b=a*[b*c*b*c*b*c*b*a*b*c*b*c*b*c*b]=
a*(a*(c*(c*b)))=a$;
\item[(15)] Follows from (14), since \\
$a*[b*c*b*c]*a*c*b*c*a=a*c*b*c*b*a*c*b*c*a$;
\item[(16)] Since $a*b*c*a*b*c*a=^{(5)}a*c*b*a*c*b$, we have:\\
$a*b*[c*b*c]*a*b*[c*b*c]*a=a*b*(b*c)*a*b*(b*c)*a=a*[(b*c)]*b*a*[(b*c)]*b=
a*[c*b*c*b]*a*[c*b*c*b]=a*b*c*b*c*a*b*c*b*c$.
\end{enumerate}
One of the oldest conjectures concerning local changes in the diagram is the Nakanishi conjecture (see \cite{Prz}).
\begin{conjecture}[Nakanishi, 1979]
Every knot is $4$-move equivalent to the trivial knot.
\end{conjecture}
Our hope was that the fourth Burnside kei of the knot can be used to 
detect a potential counterexample to this conjecture. 
However, the following theorem suggests this is not likely to be the case.
\begin{theorem}
Every algebraically connected quotient of $\bar Q(3,4)$ is a trivial, 
one element quandle.
\end{theorem}
\begin{proof}
Let $\tilde{Q}$ be the algebraically connected quotient of $\bar Q(3,4)$, \\
$f\colon \bar Q(3,4) \to \tilde{Q}$ be the quotient homomorphism and $S_1$, $S_2$,
$S_3$ denote the algebraically connected components of $\bar Q(3,4)$ .\\
We claim that $\tilde{Q}$ is contained in each image $f(S_i)$, for $i=1,\,2,\,3$.\\
By the way of contradiction, let us assume that there exists $x\in\tilde{Q}$ and
\mbox{$f^{-1}(x)\,\cap\, S_j=\emptyset$,}
for some $j$. Let $a$ be any element of $S_j$ and $y=f(a)$. Then from the algebraic connectivity of $\tilde Q$ follows that
$x=y*x_1*\ldots *x_k$, for some elements $x_1,\ldots,\,x_k\in\tilde{Q}$.
Now we choose arbitrary elements $z_i$ from the preimages $f^{-1}(x_i)$, 
$i=1,\ldots,\, k$. Let $z=a*z_1*\ldots *z_k$. Then $z\in S_j$ and
$f(z)=f(a)*f(z_1)*\ldots*f(z_k)=y*x_1*\ldots*x_k=x$, 
which contradicts the original assumption 
that $f^{-1}(x)\,\cap\, S_j=\emptyset$.\\
Each subquandle $S_i$, when considered as a quandle itself, 
has eight $4$-element components (orbits), $T_1,\ldots,T_8$ and 
just as before, we can prove that $\tilde{Q}$ is the image of each 
$T_i$ (and can have at most 4 elements). But finally we can use the fact 
that every such $T_i$ is a trivial quandle ($x*y=x$, for any $x$, 
$y\in T_i$) and $\tilde{Q}$ must be an image of just one element.
\end{proof}
Since knot quandles are algebraically connected, we have 
the following result.
\begin{corollary}
Let $K$ be a knot such that the minimal number of generators 
of its fundamental quandle is
$\leq 3$. Then its fourth Burnside quandle, $\bar Q_4(K)$, has only 
one element.
\end{corollary}
For example, $\bar Q_4(K)$ will not detect a potential counterexample 
to the Nakanishi $4$-move conjecture among $3$-bridge knots.\\
It seems to be plausible that the order of $\bar Q_4(K)$ is $1$, 
for any knot $K$.
 
\begin{center}
\begin{tabular}{l  @{\hspace{2.5 cm}} l}
                            &                         \\
Maciej Niebrzydowski        &     J\'ozef~H. Przytycki\\
e-mail: niebrz@gwu.edu      &     e-mail: przytyck@gwu.edu\\
                            &  
\end{tabular}
\end{center}
\noindent \textsc{Dept. of Mathematics, Old Main Bldg., 1922 F St. NW\\
The George Washington University, Washington, DC 20052}

\begin{thebibliography}{99}
\bibitem[Bu]{Bu} W.~Burnside, On an Unsettled Question in the Theory of
Discontinuous Groups. {\it Quart. J. Pure Appl. Math.} 33, 1902, 230-238.
\bibitem[C-E-S]{C-E-S} J.~S. Carter, M. Elhamdadi, M. Saito,
Twisted quandle homology theory and cocycle knot invariants,
\textit{Algebraic \& Geometric Topology}, Volume 2, 2002, 95-135.
\bibitem [D-I-P]{D-I-P}
M.~K. D{\c a}bkowski, M. Ishiwata, J.~H. Przytycki,
Rational moves and tangle embeddings: (2,2)-moves as a case study,
{\it Proceedings of the Conference Topology of Knot VII}
(held at TWCU, December 23-26, 2004), February, 2005, 37-46, in Japanese;\\
e-print (in English): http://front.math.ucdavis.edu/math.GT/0501539

\bibitem[F-R]{F-R} R. Fenn, C. Rourke, Racks and links in codimension two,
\textit{J. of Knot Theory and its Ramifications}, Volume 1, No 4, 343-406,
World Scientific 1992. 
\bibitem[Joy]{Joy} D.Joyce, A classifying invariant of knots: the knot
quandle, \textit{J. Pure Appl. Alg.}, 23, 1982, 37-65.
\bibitem [Kam]{Kam} S. Kamada, Knot invariants derived from quandles and racks,
\textit{Geometry \& Topology Monographs}, Volume 4, 2002-4, 103-117;\\
http://front.math.ucdavis.edu/math.GT/0211096
\bibitem[Oht]{Oht} T. Ohtsuki, personal communication, January 2005.
\bibitem[Prz]{Prz} J.~H. Przytycki, $t_k$-moves on links,
{\it Braids}, ed. J.~S. Birman and A. Libgober, Contemporary Math., 
Volume 78, 1988, 615-656.
\bibitem[Tak]{Tak}
M. Takasaki, Abstraction of symmetric transformation, (in Japanese) 
{\it Tohoku Math. J.}, 49, 1942/3, 145-207.
\bibitem[VL-1]{VL-1} M. Vaughan-Lee, {\it The restricted Burnside problem};
Second edition. London Mathematical Society Monographs. New Series, 8. The
Clarendon Press, Oxford University Press, New York, 1993. xiv+256 pp.
\bibitem[Win]{Win}
S.~K. Winker, Quandles, knot invariants and the n-fold branched cover,
PhD thesis, University of Illinois at Chicago, 1984.
\end{thebibliography}
\end{document}